\documentclass[aop,MSNbibl,seceqn,citesort,dvips]{arximspdf}

\doi{10.1214/11-AOP678}
\volume{40}
\issue{5}
\pubyear{2012}
\firstpage{2106}
\lastpage{2130}

\makeatletter

\newcommand{\eps}{\varepsilon}
\newcommand{\Z}{{\mathbb Z}}
\newcommand{\N}{{\mathbb N}}
\newcommand{\M}{{\mathcal M}}
\newcommand{\D}{{\mathcal D}}
\newcommand{\MM}{{\mathbb M}}
\newcommand{\G}{{\mathcal G}}
\newcommand{\R}{{\mathbb R}}
\newcommand{\RR}{{\mathcal R}}
\renewcommand{\phi}{\varphi}

\newcommand{\Sph}{\mathbb{S}}

\newcommand{\PV}{\mathbf{P}}
\newcommand{\IE}{\mathbb{E}}
\newcommand{\IP}{\mathbb{P}}

\newcommand{\F}{\mathcal{F}}

\newtheorem{theo}{Theorem}[section]
\newtheorem{lmm}[theo]{Lemma}
\newtheorem{prop}[theo]{Proposition}
\newtheorem{cor}[theo]{Corollary}

\newproclaim{rmk}{Remark}[section]

\newproclaim{Condition}{Condition}
\newproclaim{Conditions}{Conditions}

\makeatother

\begin{document}
\begin{frontmatter}

\title{On a general many-dimensional excited random~walk}
\runtitle{General many-dimensional excited random walk}

\begin{aug}
\author[A]{\fnms{Mikhail} \snm{Menshikov}\thanksref{t1}\ead[label=e1]{Mikhail.Menshikov@durham.ac.uk}},
\author[B]{\fnms{Serguei} \snm{Popov}\corref{}\thanksref{t2}\ead[label=e2]{popov@ime.unicamp.br}},
\author[C]{\fnms{Alejandro F.} \snm{Ram\'\i rez}\thanksref{t3}\ead[label=e3]{aramirez@mat.puc.cl}}\\ and
\author[B]{\fnms{Marina} \snm{Vachkovskaia}\thanksref{t2}\ead[label=e4]{marinav@ime.unicamp.br}}
\runauthor{Menshikov, Popov, Ram\'\i rez and Vachkovskaia}
\affiliation{University of Durham, University of Campinas, Pontificia
Universidad Cat\'olica de Chile
and University of Campinas}
\address[A]{M. Menshikov\\
Department of Mathematical Sciences\\
University of Durham\\
South Road, Durham DH1 3LE\\
United Kingdom\\
\printead{e1}} 
\address[B]{S. Popov\\
M. Vachkovskaia\\
Department of Statistics\\
Institute of Mathematics, Statistics\\
\quad and Scientific Computation\\
University of Campinas---UNICAMP\\
Rua S\'ergio Buarque de Holanda 651\\
13083-859, Campinas SP\\
Brazil\\
\printead{e2}\\
\hphantom{E-mail: }\printead*{e4}}
\address[C]{A. F. Ram\'{\i}rez\\
Facultad de Matem\'aticas\\
Pontificia Universidad Cat\'olica de Chile\\
Vicu\~na Mackenna 4860\\ Macul
Santiago 6904441\\
Chile \\
\printead{e3}}
\end{aug}

\thankstext{t1}{Supported in part by Mecesup
0711 during his stay at the Department of Mathematics,
Pontificia Universidad Cat\'olica de Chile, where part of this work
was done.}

\thankstext{t2}{Supported by CNPq Grants 300886/2008--0,
301455/2009--0, 472431/2009--9
and by FAPESP thematic Grant 09/52379--8.}

\thankstext{t3}{Supported in part by
Fondo Nacional de Desarrollo Cient\'\i fico y Tecnol\'ogico
Grant 1100298.}

\received{\smonth{1} \syear{2010}}
\revised{\smonth{10} \syear{2010}}

%
\begin{abstract}
In this paper we study a substantial generalization of the model of
excited random walk introduced in [\textit{Electron. Commun. Probab.}
\textbf{8} (2003) 86--92] by Benjamini and Wilson. We
consider a discrete-time stochastic process $(X_n,n=0,1,2,\ldots)$
taking values on $\Z^d$, \mbox{$d\geq2$}, described as follows: when the
particle visits a site for the first time, it has a uniformly-positive
drift in a given direction $\ell$; when the particle is at a site which
was already visited before, it has zero drift. Assuming uniform
ellipticity and that the jumps of the process are
uniformly bounded, we prove that the process is ballistic in the
direction $\ell$ so that $\liminf_{n\to\infty}\frac{X_n\cdot
\ell}{n}>0$. A~key ingredient in the proof of this result is an
estimate on the probability that the process visits less than
$n^{{1/2}+\alpha}$ distinct sites by time $n$, where $\alpha$ is
some positive number depending on the parameters of the model. This
approach completely avoids the use of tan points and coupling methods
specific to the excited random walk. Furthermore, we apply this
technique to prove that the excited random walk in an i.i.d. random
environment satisfies a ballistic law of large numbers and a central
limit theorem.
\end{abstract}

%

\end{frontmatter}

\section{Introduction and results}
\label{sintro}
Let $p\in(1/2,1]$. Consider two discrete time simple random walks
on the hyper-cubic lattice $\Z^d$, $d\geq2$: a symmetric
random walk $(Y_n,n\ge0)$ and a random walk $(Z_n,n\ge0)$ which
jumps to the right with probability $p/d$, to the left with
probability $(1-p)/d$
and to the other nearest-neighbor sites with probability $1/(2d)$.
The \textit{excited} or \textit{cookie random walk} with bias
parameter $p$
on $\Z^d$ is a self-interacting random walk $(X_n,n\ge0)$
starting from $0$ and defined as follows: if at time $n$
the walk is at a site $x$ which it visited at some time $k$ such
that $k<n$,
it jumps according to the transition probabilities of the symmetric
random walk $(Y_n)$, so that it jumps
with probability $1/(2d)$ to the nearest-neighbor sites of $x$;
if at time $n$ the process visits a site $x$ for the first time,
it jumps according to the transition probabilities of the walk
$(Z_n)$, so that it
jumps to the right with probability $p/d$, to the left with
probability $(1-p)/d$ and to the other nearest-neighbor sites
of $x$ with probability $1/(2d)$ (eating one \textit{cookie}
at site~$x$). We will call the walks $(Y_n)$
and $(Z_n)$ defining the excited random walk, the
\textit{underlying processes} of the excited random walk.

The excited random walk was introduced in 2003 by Benjamini and
Wilson in~\cite{BW03}. They proved that for dimensions $d\ge2$
it is transient to the right, meaning that a.s.
%
\begin{equation}
\label{transience}
\lim_{n\to\infty}X_n\cdot e_1=\infty,
\end{equation}
where $\{e_i\dvtx1\le i\le d\}$ denote the canonical generators of the
additive group~$\Z^d$.
Furthermore, they showed that in dimensions $d\ge4$ the excited
random walk is ballistic to the right so that a.s.
%
\begin{equation}
\label{ballisticity}
\liminf_{n\to\infty}\frac{X_n\cdot e_1}{n}>0.
\end{equation}
In~\cite{K03} and~\cite{K05}, Kozma extended the
above ballisticity result to dimensions $d=3$ and $d=2$.
Standard methods based on regeneration times
can be used to deduce from (\ref{ballisticity})
that a law of large numbers with deterministic speed~$v$ is
satisfied, so that a.s.
\[
\lim_{n\to\infty}\frac{X_n}{n}=v,
\]
where $v\cdot e_1>0$ and $v\cdot e_j=0$ for $2\le j\le n$.
In~\cite{BR07}, B\'erard and Ram\'\i{}rez gave an alternative proof
of ballisticity and proved
that a central limit theorem is satisfied in dimensions $d\ge2$,
so that
%
\begin{equation}
\label{clt}
\eps^{1/2}(X_{\eps^{-1}n}-\eps^{-1}nv)
\end{equation}
converges in law as $\eps\to0$ to a Brownian motion with
variance $\sigma^2>0$. A variant of the excited random walk, called
the \textit{multi-excited random walk}, was introduced by Zerner
in~\cite{Z05},
where the walk has the possibility of consuming more
than one cookie per site, and hence the process exhibits a
nontrivial
behavior even in dimension $d=1$. Several papers have been written
where the
transient and ballisticity properties of this model are studied in
random
and deterministic environments, mainly in dimension $d=1$
(see, e.g.,~\cite{KZ08,KM09,BS09}).
Nevertheless, with the exception of~\cite{KRS09} and~\cite{Z06}, a
very natural issue has not been so far addressed:
what happens if the underlying processes are no longer
nearest-neighbor spatially homogeneous random
walks? For example, if $(Y_n)$ and $(Z_n)$ are random walks on
$\Z^d$, which are not spatially homogeneous and do not perform
nearest-neighbor jumps, $(Y_n)$ has zero drift, and $(Z_n)$ has a~drift
to the right, is the corresponding excited random walk transient
to the right?
It is reasonable to wonder under which conditions the corresponding
excited random walk would still be transient to the right as
in (\ref{transience}),
ballistic as in (\ref{ballisticity}), or satisfy a
central limit theorem as in (\ref{clt}).
In this paper, we study a generalization of the excited random walk
on $\mathbb Z^d$, $d\geq2$, where the underlying processes of the
model are not necessarily homogeneous random walks, and not even
Markovian.

From our point of view, part of the reason why these issues have not
been considered is related to the techniques so far developed to
study the many-dimensional excited random walk. Indeed,
the proof of transience to the right of~\cite{BW03}, the
law of large numbers and the central
limit theorem of~\cite{BR07} in dimensions $d\ge2$,
rest on the following two key ingredients:
(i) the excited random walk
can be coupled to the underlying simple symmetric random walk
$(Y_n)$, in such a way that $(X_n-Y_n)\cdot e_1$ is nondecreasing
in $n$ and $(X_n-Y_n)\cdot e_j=0$ for every $2\le j\le d$;
(ii) it is possible to get a lower bound for the cardinality of
the range at time $n$ of the excited random walk in terms of the
\textit{tan points} of the coupled simple symmetric random walk,
which exploit reversibility properties of the symmetric random walk
(see~\cite{BW03} and~\cite{BMS02}).
A tan point of the random walk $(Y_n)$
in dimension $d=2$
is defined as any site $x\in\Z^2$ with the property that
the ray $\{x+ke_1\dvtx k\ge0\}$ is visited
by the random walk for the first time
at site $x$. As explained in~\cite{BR07}, it can be shown using
the ideas of~\cite{BW03} that for every $\eps$
the number of tan points visited by $(Y_m)$ at time $n$ is larger
than $n^{3/4-\eps}$ with a probability that decays faster than any
polynomial in $n$.
Notwithstanding, these methods break down when the underlying
processes are even slightly modified.
On one hand, the coupling between the excited random walk
and the simple random walk does not work in general. Furthermore,
the estimation of the number of tan points at a given
time is a very specific argument which works only
for simple symmetric random walks.
Hence, to study excited random walks defined in
terms of more general underlying processes,
more powerful
methods have to be developed: a fundamental
ingredient of this paper is the introduction of a new technique to
estimate the range of general versions of excited random walks, which
completely avoids the use of tan points.

We define a \textit{generalized excited random
walk} which will correspond to a process driven
by an underlying process $(Y_n)$ which is a ($d$-dimensional)
martingale and
a process $(Z_n)$ which satisfies minimal requirements, including
the presence of a drift. We develop a machinery
which avoids the use of
tan points and\vspace*{1pt} of coupling, proving that this generalized excited
random walk is ballistic.
Let \mbox{$\|\cdot\|$} be the $L^2$-norm in $\Z^d$ or $\R^d$, $d\geq2$;
also, we define $\Sph^{d-1}=\{x\in\R^d \dvtx \|x\|=1\}$ to be the
unit sphere in $\R^d$.
Consider a~$\Z^d$-valued stochastic process
$X=(X_n, n=0,1,2,\ldots)$ adapted to a filtration
$\F=(\F_n, n=0,1,2,\ldots)$.
Unless otherwise stated, we suppose that $X_0=0$.
Denote by $\IP$ the law of $X$ and by $\IE$ the corresponding
expectation.
As mentioned before, the processes we are considering are
known also as \textit{cookie random walks}.
This terminology is also useful to us, because in the sequel we will
need to consider situations when the particle gets the first visit
push not in all the sites, but only in the sites of some fixed subset
of $\Z^d$. In this case, we say that the initial configuration of
cookies (or the initial cookie environment) is such that they are
only in this set.

Throughout the paper we suppose that
the jumps of the process are uniformly bounded, that is,
the following condition holds for the process $X$:

\renewcommand{\theCondition}{B}
\begin{Condition}\label{ConditionB}
There exists a constant $K>0$ such that
${\sup_{n\ge0}}\|X_{n+1}-X_n\|\leq K$ a.s.
\end{Condition}


Next, consider the following condition:
\renewcommand{\theCondition}{C$^+$}
\begin{Condition}\label{ConditionC+}
Let $\ell\in\Sph^{d-1}$. We say that Condition~\ref{ConditionC+} is
satisfied with respect to $\ell$ if there exist a $\lambda> 0$
such that
\[
\IE(X_{n+1}-X_n\mid\F_n)=0 \qquad\mbox{on }
\{\mbox{there exists }k<n \mbox{ such that }X_k=X_n\}
\]
and
\[
\IE(X_{n+1}-X_n\mid\F_n)\cdot\ell\geq\lambda\qquad\mbox{on }
\{X_k\neq X_n \mbox{ for all }k<n\}.
\]
\end{Condition}

The meaning of Condition~\ref{ConditionC+}
is that, when the process $X$
visits a site for the first time, it has drift in the
direction $\ell$, whereas if it comes to an already visited
site, it has zero drift behaving like a martingale.


Also, we formulate:
\renewcommand{\theCondition}{E}
\begin{Condition}\label{ConditionE}
Let $\ell\in\Sph^{d-1}$. We say that Condition~\ref{ConditionE} is satisfied with
respect to $\ell$ if there exist $h,r>0$ such that for all $n$
%
\begin{equation}
\label{snos}
\IP[(X_{n+1}-X_n)\cdot\ell>r\mid\F_n] \geq h
\end{equation}
and for all $\ell'$ with $\|\ell'\|=1$, on
$\{\IE(X_{n+1}-X_n\mid\F_n)=0\}$
%
\begin{equation}
\label{snoss}
\IP[(X_{n+1}-X_n)\cdot\ell'>r\mid\F_n] \geq h.
\end{equation}
\end{Condition}

Condition~\ref{ConditionE} is a kind of uniform
ellipticity assumption which states that the process can always
advance in the direction $\ell$ by a uniformly positive amount with
a uniformly positive probability, and also, when the local drift is
equal to zero, the process can do so in any direction.
In fact, one may easily verify that if the Conditions~\ref{ConditionB} and
\ref{ConditionC+}
are satisfied, then automatically (\ref{snos}) holds for some
positive $h,r$.
However, we still formulate Condition~\ref{ConditionE} this way because in the
sequel we will need also to consider processes where the
first visit push is not necessarily uniformly positive.
Now, given $\ell\in\Sph^{d-1}$, any stochastic process $X$ adapted
to a filtration $\F$, which
satisfies Condition~\ref{ConditionB} and Conditions~\ref{ConditionC+} and~\ref{ConditionE} with respect
to~$\ell$, will be called a
\textit{generalized excited random walk} in direction $\ell$.
Note that the standard
excited random walk with its natural filtration
is a generalized random walk in direction $e_1$.
The first result of this paper is that any generalized
excited random walk is ballistic.
\begin{theo}
\label{tballisticity}
Let $d\geq2$ and $\ell\in\Sph^{d-1}$.
Assume that $X$ is a generalized excited
random walk in direction $\ell$. Then,
there exists $v=v(d,K,h,r,\lambda)>0$ such that
%
\begin{equation}
\label{eqballisticity}
\liminf_{n\to\infty} \frac{X_n\cdot\ell}{n}\geq v \qquad\mbox{a.s.}
\end{equation}
\end{theo}

The proof of Theorem~\ref{tballisticity} rests on two
key ingredients: a general result which
says that a $d$-dimensional martingale
satisfying Condition~\ref{ConditionB} and a~condition analogous to Condition~\ref{ConditionE}
should typically visit much more than~$t^{1/2}$ distinct sites by
time $t$; then, the same kind of result is also obtained
for any generalized excited random walk with an arbitrary initial
configuration of cookies. It is worth noting that
the approach of this paper could be applied also to models with
other rules of assigning the drift to the particle. For instance,
one can consider a model mentioned in~\cite{Z06}: the random walk
receives a push in the direction $\ell$ not at the
first visit, but at the $k$th visit to the site, where $k>1$.
To study such a model, one needs to prove that the set of sites
\textit{visited at least $k$ times} should be sufficiently large.
However, it is not
difficult to do so using (some suitable) uniform ellipticity
condition: if the particle visits a site, then in the next few
instants of time this site will be visited $k-1$ times more with a
uniformly positive probability.

In this paper, we also consider an excited random
walk in an i.i.d. random environment in $\Z^d$, $d\geq2$,
proving a law of large numbers and a central limit theorem for it.
Let $\mathcal P$ be the set of probability measures on
$\{\pm e_i, 1\le i\le d\}$. Let $\M=\mathcal P^{\N}$ and $\Omega=
\M^{\mathbb Z^d}$. An element $\omega=\{\omega(x),x\in\mathbb
Z^d\}\in\Omega$ is called
an \textit{environment}. Here, for each $x\in\mathbb Z^d$,
$\omega(x)=\{\omega_n(x),n\ge0\}\in\M$ and
$\omega_n(x)=\{\omega_n(x,e)\in
\mathcal P\}$. Let $\PV$ be a probability measure defined
on $\Omega$ under which the random variables $\{\omega(x),x\in
\mathbb
Z^d\}$
are i.i.d. Let us stress that we do not assume any
independence of the random variables $\omega_n(x)$, $n\ge0$, for a
fixed $x$. We assume that $\PV$ is \textit{uniformly elliptic} so that
there exists a constant $\kappa>0$ such that for
every $n\ge0$ and $1\le i\le d$ one has that
$\PV[\omega_n(0,\pm e_i)\ge\kappa]=1$. Furthermore, we assume that
$\PV$ is \textit{uniformly excited} in the direction $\ell\in\mathbb
S^{d-1}$,
so that there exists a $\lambda>0$ such that
\[
\PV\Biggl[\sum_{i=1}^d\bigl(e_i \omega_0 (0,e_i)-e_i\omega_0
(0,-e_i)\bigr)\cdot\ell\ge\lambda\Biggr]=1,
\]
and such that for every $j\ge1$ we have that
\[
\PV\Biggl[\sum_{i=1}^d\bigl(e_i \omega_j (0,e_i)-e_i\omega_j
(0,-e_i)\bigr)=0\Biggr]=1.
\]
We now define the \textit{excited random walk in random
environment} (ERWRE) as
the process defined for each $n\ge0$, $x\in\mathbb Z^d$ and
$e\in\mathbb Z^d$ with $|e|=1$, through the transition probabilities,
\[
P_\omega\Biggl[X_{k+1}=x+e \Bigm| X_k=x,
\sum_{j=0}^{k-1}{\mathbf{1}}_{\{X_j=x\}}=n\Biggr] = \omega_n(x,e).
\]
In other words, whenever the process visits a site $x$ for the first
time, it has a~mean drift in the direction $\ell$ which is larger
than $\lambda$; whenever it visits a site~$x$ which it already
visited before, its mean drift is $0$. For each environment~$\omega$
and $x\in\mathbb Z^d$, we call $P_{\omega,x}$ the law of such a
process starting from $x$. We define the \textit{annealed} or
\textit{averaged} law of
the excited random walk in random environment as $P_x=\int
P_{\omega,x}\,d\PV$, in opposition to $P_{\omega,x}$ which is called
the \textit{quenched law}. The second result of this paper is the
following theorem.
\begin{theo}
\label{trwre}
Consider the excited random walk in random
environment, uniformly
excited in the direction $\ell\in\mathbb S^{d-1}$.
Then, the following are satisfied:
\begin{longlist}[(ii)]
\item[(i)] \textup{(Law of large numbers)}. There exists $v$ such
that
$v\cdot\ell>0$ and
\[
\frac{X_n}{n}\to v, \qquad\mbox{$P_0$-a.s.}
\]
\item[(ii)] \textup{(Central limit theorem)}. There exists a
nondegenerate matrix $A$ such that
\[
\eps^{1/2}\bigl(X_{\lfloor n\eps^{-1}\rfloor}-n\eps^{-1}v\bigr)
\]
converges as $\eps\to0$ in $P_0$-law to a $d$-dimensional Brownian
motion with covariance matrix $A$.
\end{longlist}
\end{theo}

In the above theorem, we considered the case of nearest-neighbor
jumps only for notational convenience; it extends to the case of
uniformly bounded jumps without difficulties (one has to assume
also that\vadjust{\goodbreak} Condition~\ref{ConditionE} holds).
Observe also that, taking $\PV$ concentrated in one point (so that
the environment in all sites is the same), we obtain LLN and CLT for
a spatially homogeneous generalized excited random walk.

Let us note that, in~\cite{Z06} transience is proved for
an excited random walk in random environment using the method of
environment viewed from the particle. As explained in~\cite{Z05},
this implies using regeneration times, a law of large numbers with
a velocity which could be possibly equal to $0$. Nevertheless,
we do not see how to use the techniques based on coupling with
a~simple symmetric random walk and tan points, to prove that
the expected value of the regeneration times is finite.

This paper is organized in the following way. First, in
Section~\ref{srenewalstr} we define a sequence of regeneration
times. The key fact about this sequence (which is essential for
proving Theorem~\ref{tballisticity}) is that the time intervals
between the regenerations behave in a nice way; see
Proposition~\ref{tail-reg}. The proof of this proposition is rather
technical and is postponed to Section~\ref{sprooftails}
(generally, in this paper we prefer to postpone the proofs of more
technical results). In Section~\ref{smain} we prove
Theorem~\ref{tballisticity} and then apply the previously developed
machinery to the excited random walk in random
environment proving Theorem~\ref{trwre}. In
Section~\ref{sdisplest} we study typical
displacement of excited random walk by time $n$; these estimates are
then used in Section~\ref{sprooftails}. It turns our that, to
understand what the typical displacement of the excited random walk
should be, one has to understand the typical behavior of the number
of different sites visited by time $n$ (i.e., the range of the
process). The key result concerning the range
(Proposition~\ref{phitpoints}) is stated in
Section~\ref{sdisplest} without proof. In
Section~\ref{sprelimfacts} we use some martingale techniques to
obtain several auxiliary facts, which then are used in
Section~\ref{sproofhitpoints} to
prove Proposition~\ref{phitpoints}.

\section{The renewal structure}
\label{srenewalstr}
The proof of Theorems~\ref{tballisticity} and~\ref{trwre},
uses classical renewal time methods. Let us note that the
ERWRE, in a fixed environment, is a~generalized excited random walk.
Hence, here we will focus on the construction of the renewal
structure for a generalized excited random walk,
following the standard approach and notation
presented in~\cite{BR07} and in the context of random walk
in random environment in~\cite{SZ99}. Due to the fact that the
generalized excited random walk is neither space homogeneous nor
Markovian,
we will need to introduce a general notation, and deviate slightly
from the construction of~\cite{BR07}.
Let $\ell\in\Sph^{d-1}$. We consider a stochastic process
$(X_n,n\ge0)$ satisfying Conditions~\ref{ConditionB},~\ref{ConditionE} and~\ref{ConditionC+},
with respect to a~filtration~$\F_n$ and a direction
$\ell\in\Sph^{d-1}$.
For each $u>0$ let
\[
T_u=\min\{k\ge1\dvtx X_k\cdot\ell\ge u\}.
\]
Define
\[
\bar D=\inf\{m\ge0\dvtx X_m\cdot\ell< X_0\cdot\ell\}.
\]
Furthermore, define two sequences
of $\F_n$-stopping times $\{S_{n}\dvtx n\ge0\}$
and $\{D_{n}\dvtx n\ge0\}$ as follows. We\vadjust{\goodbreak} let $S_{0}=0$,
$R_{0}=X_0\cdot\ell$ and $D_{0}=0$.
Next, define by induction in $k\ge0$
\begin{eqnarray*}
S_{k+1}&=&T_{R_{k}+1},\\
D_{k+1}&=&\bar D\circ\theta_{S_{k+1}}+S_{k+1},\\
R_{k+1}&=&\sup\{X_i\cdot\ell\dvtx0\le i\le D_{k+1}\},
\end{eqnarray*}
where $\theta$ is the canonical shift on the space of trajectories.
Let
\[
\kappa=\inf\{n\ge0\dvtx S_{n}<\infty, D_{n}=\infty\}
\]
with the convention that $\kappa=\infty$ when
$\{n\dvtx S_{n}<\infty, D_{n}=\infty\}=\varnothing$.
We define the first regeneration time as
\[
\tau_1=S_{\kappa}.
\]
We then define by induction on $n\ge1$, the sequence
of regeneration times $\tau_1,\tau_2,\ldots$ as follows:
\[
\tau_{n+1}=\tau_n+\tau_1(X_{\tau_n+\cdot})
\]
setting $\tau_{n+1}=\infty$ when $\tau_n=\infty$.

Next, we define $D^{(0)}_i=D_{i}$ and
$S^{(0)}_i=S_{i}$ and
for each $k\ge1$ two sequences $D^{(k)}_i$ and
$S^{(k)}_i$ corresponding to the regeneration time $\tau_{k+1}$,
analogously to the definition of the sequences of stopping
times $D_{i}$ and $S_{i}$ related to $\tau_1$. For example,
$S^{(1)}_0=\tau_1$, $R^{(1)}_0=X_{\tau_1}\cdot\ell, D^{(1)}_0=0$
and define by induction in $i\ge0$,
\begin{eqnarray*}
S^{(1)}_{i+1}&=&T_{R^{(1)}_i+1},\\
D^{(1)}_{i+1}&=&\bar D\circ\theta_{S^{(1)}_{i+1}}+S^{(1)}_{i+1},\\
R^{(1)}_{i+1}&=&\sup\bigl\{X_i\cdot\ell\dvtx0\le i\le D^{(1)}_{i+1}\bigr\}.
\end{eqnarray*}
As opposed to the situation which occurs for the standard
excited random walk (see~\cite{BR07}), here the sequence
of regeneration times is not necessarily i.i.d.
For each $k\ge1$ and $j\ge0$ such that
$S^{(k)}_j<\infty$, we need to introduce the $\sigma$-algebra of
events
up to time $S^{(k)}_j$. We define $\G^{(k)}_j$ as the
smallest $\sigma$-algebra
containing all sets of the form
$\{\tau_1\le n_1\}\cap\cdots\cap\{\tau_k\le n_k\}\cap A$, where
$n_1<n_2<\cdots< n_k$ are integers and $A\in\F_{n_k+S_j^{(k)}}$.
In Section~\ref{srange}, we will prove the
following proposition, which will play a key role in the
proof of Theorem~\ref{tballisticity}.
\begin{prop}
\label{tail-reg}
Consider a generalized random walk excited
in the direction $\ell$, and let $(\tau_k,k\ge1)$ be the
associated sequence of regeneration times.
Then,
there exist $C',\alpha'>0$ such that for every $n\ge1$,
\[
\sup_{k\ge1}\IP\bigl[\tau_{k+1}-\tau_k>n \mid\G^{(k)}_0\bigr] \le
C'e^{-n^{\alpha'}} \qquad\mbox{a.s.}
\]
In particular, for every $k\ge1$ we have that $\tau_k<\infty$ a.s.\vadjust{\goodbreak}
\end{prop}

Now, let us prove the following result which will be useful in the
sequel.
Throughout, we denote by $\varpi$ an element of the space
of trajectories $(\Z^d)^{\N}$. Furthermore, we define for $n\ge1$,
\[
\bar D_n:=\inf\{m\ge n\dvtx X_m\cdot\ell<X_n\cdot\ell\}.
\]

\begin{prop}
\label{conditioning}
Let $A$ be a Borel subset of $(\Z^d)^{\N}$.
Then, the following statements are satisfied:
\begin{longlist}[(ii)]
\item[(i)] For every $k\ge1$,
\[
\IP\bigl[X_{\tau_k+\cdot}\in A \mid\G^{(k)}_0\bigr]=\sum_{n=1}^\infty
{\mathbf{1}}_{\{\tau_k=n\}}(\varpi)\IP[X_{n+\cdot}\in A \mid\bar
D_n=\infty
,\F_n].
\]
\item[(ii)] For every $k,j\ge1$,
\begin{eqnarray*}
&&\IP\bigl[X_{S^{(k)}_j+\cdot}\in A \mid\G^{(k)}_j\bigr]\\
&&\qquad=\sum_{n=1}^\infty
\sum_{m=1}^\infty{\mathbf{1}}_{\{\tau_k=n\}} (\varpi)
{\mathbf{1}}_{\{S_j^{(k)}=n+m\}}(\varpi)\IP[X_{n+m+\cdot}\in A \mid
\bar
D_n=\infty,
\F_{n+m}].
\end{eqnarray*}
\end{longlist}
\end{prop}
\begin{pf}
Since the proof of part (i) is simpler than that
of part (ii), we will omit it. For part (ii), we only
consider the case $k=1$ and $j=1$,
being the case $k=1$ and $j>1$ similar.
Using the fact that for every natural $n$ and $k>1$, the event
$\{\bar D_n=\infty\}\cap\{n<\tau_k\}$ is $\mathcal G_0^{(k)}$
measurable,
the cases when $k>1$ can be proved in a similar way.
For each $n,m\ge1$, define $\mathcal S_{n,m}$ as the set of
trajectories $\{x_0,\ldots,x_{n-1},x_n,x_{n+1},\ldots,x_{n+m}\}$
satisfying the following properties:
\begin{longlist}[(a)]
\item[(a)]
\begin{eqnarray*}
\\[-30pt]
x_n\cdot\ell>\sup_{0\le l\le n-1}x_l\cdot\ell;
\end{eqnarray*}
\item[(b)] for each $l$ such that $0\le l\le n-1$ one has
\[
\min_{T_{x_l\cdot\ell}< i\le n-1}x_i\cdot\ell< x_l\cdot\ell;
\]
\item[(c)] one has that
\[
x_{n+m}\cdot\ell\ge x_n\cdot\ell+1>\sup_{0\le l\le
n+m-1}x_l\cdot\ell.
\]
\end{longlist}
These three conditions define the trajectories in
$\mathcal S_{n,m}$ as
those for which if
\mbox{${\bar D}_n=\infty$}, then
$S^{(1)}_1=n+m$ (and $\tau_1=n$).
We will use the notation $s_{n,m}$ for an element of
$\mathcal S_{n,m}$. Furthermore, given a trajectory $\varpi\in
(\Z^d)^{\N}$,
we will denote by $\varpi_n$ its projection to the first $n$
coordinates.
Now note that $\mathcal G_1^{(1)}$ is
generated by the disjoint collection of sets of the form
\[
\{\varpi_{n+m}=s_{n,m}\}\cap\{\tau_1=n\},\vadjust{\goodbreak}
\]
where $n$ and $m$ vary over the naturals and
$s_{n+m}\in\mathcal S_{n,m}$.
Hence
\begin{eqnarray*}
&&\IP\bigl[X_{S^{(1)}_1+\cdot}\in A \mid\mathcal G_1^{(1)}\bigr]
\\
&&\qquad=
\sum_{n=1}^\infty{\mathbf{1}}_{\{\tau_1=n\}}
\sum_{m=1}^\infty\sum_{s_{n+m}\in\mathcal S_{n,m}}
{\mathbf{1}}_{\{\varpi_{n+m}=s_{n+m}\}}\\
&&\qquad\quad\hphantom{\sum_{n=1}^\infty{\mathbf{1}}_{\{\tau_1=n\}}
\sum_{m=1}^\infty\sum_{s_{n+m}\in\mathcal S_{n,m}}}
{}\times\IP[X_{n+m+\cdot}\in
A\mid\tau_1=n,
\varpi_{n+m}=s_{n+m}].
\end{eqnarray*}
On the other hand, for each $n$ and $m$ we have that when
$s_{n+m}\in\mathcal S_{n,m}$
\[
\{\varpi_{n+m}=s_{n+m}\}\cap\{\tau_1=n\}=\{\varpi_{n+m}=s_{n,m},
\bar D_n=\infty\}.
\]
Therefore, since
${\mathbf{1}}_{\{\tau_1=n\}}(\varpi){\mathbf{1}}_{\{S_1^{(1)}=n+m\}
}(\varpi)
{\mathbf{1}}_{\{\varpi_{n+m}=s_{n+m}\}}=0$
whenever $s_{n+m}\notin\mathcal S_{n,m}$, we see that
\begin{eqnarray*}
&&\IP\bigl[X_{S_1^{(1)}+\cdot}\in A\mid\mathcal G_1^{(1)}\bigr]\\
&&\qquad =\sum_{n=1}^\infty\sum_{m=1}^\infty \sum_{s_{n,m}\in\mathcal
S_{n,m}}{\mathbf{1}}_{\{\tau_1=n\}}(\varpi)
{\mathbf{1}}_{\{S_1^{(1)}=n+m\}}(\varpi)\\
&&\quad\qquad\hphantom{\sum_{n=1}^\infty\sum_{m=1}^\infty \sum_{s_{n,m}\in\mathcal
S_{n,m}}}
{}\times\IP[X_{n+m+\cdot}\in A\mid \bar
D_n=\infty,\varpi_{n+m}
=s_{n,m}]\\
&&\qquad = \sum_{n=1}^\infty\sum_{m=1}^\infty{\mathbf{1}}_{\{\tau_1=n\}
}(\varpi) {\mathbf{1}}_{\{S_1^{(1)}=n+m\}}(\varpi)\\
&&\qquad\quad\hphantom{\sum_{n=1}^\infty\sum_{m=1}^\infty}
{}\times\IP[ X_{n+m+\cdot}\in
A\mid \bar D_n=\infty, \mathcal F_{n+m}].
\end{eqnarray*}
\upqed\end{pf}

\section{Proof of the main results}
\label{smain}
In this section we will prove Theorems~\ref{tballisticity}
and~\ref{trwre}.

\subsection{\texorpdfstring{Proof of Theorem \protect\ref{tballisticity}}{Proof of Theorem 1.1}}
\label{sproofmain}
To prove Theorem~\ref{tballisticity}, we first prove the following
lemma.
\begin{lmm}
\label{ltaun}
Consider a generalized excited random walk in the
direction~$\ell$. Let $(\tau_k,k\ge1)$ be the associated
regeneration times. Then, there is a~constant $C>0$ such that a.s.
\[
\limsup_{n\to\infty} \frac{\tau_n}{n}<C.
\]
\end{lmm}

\begin{pf}
Let $C'=\sup_{k\ge1}\IE(\tau_{k+1}-\tau_k \mid\G^{(k)}_0)$;
by Proposition~\ref{tail-reg}, we know that
$C'<\infty$. Let $\tau_0=0$. Now, consider the
process $M_n=\sum_{k=0}^{n-1}(\tau_{k+1}-\tau_k
-\IE(\tau_{k+1}-\tau_k \mid\G^{(k)}_0))$,
for $n\ge1$, which is a
martingale with respect to the filtration $(\mathcal G^{(n)}_0,n\ge
1)$.
We then have for $C>C'$ that
\[
\IP[\tau_n>nC] \le\IP[M_n>n(C-C')]
\le\frac{\IE[M_n^4]}{n^4 (C-C')^4}.\vadjust{\goodbreak}
\]
But, using the fact that $M_n$ is a martingale and
Proposition~\ref{tail-reg}, we see that there is
a constant $C_1>0$ such that $\IE[M_n^4 ]<C_1n^2$.
Hence,
\[
\IP[\tau_n>nC] \le\frac{C_2}{n^2}
\]
for some constant $C_2>0$,
which, by Borel--Cantelli, proves the lemma.
\end{pf}

Let us now see how to deduce Theorem~\ref{tballisticity}
from Lemma~\ref{ltaun}. By definition of the regeneration
times, note that
%
\begin{equation}
\label{xinf}
\liminf_{n\to\infty}\frac{X_{\tau_n}\cdot\ell}{n}\ge1.
\end{equation}
For each $k\ge0$, define $n_k=\sup\{n\ge0\dvtx\tau_n\le k\}$. Since
for each $n$ we have $\tau_n\ge n$,
it follows that $n_k<\infty$ a.s.
Note also that $\lim_{k\to\infty}n_k=\infty$.
Also, by definition of $\tau_k$ and $n_k$ we have
$X_k\cdot\ell\ge X_{\tau_{n_k}}\cdot\ell$, so
%
\begin{eqnarray}
\label{nk}
\liminf_{k\to\infty}\frac{X_k\cdot\ell}{k}
&=&
\liminf_{k\to\infty}\frac{X_{\tau_{n_k}}\cdot\ell}{k}
=
\liminf_{k\to\infty}\frac{n_k}{k}\frac{X_{\tau_{n_k}}\cdot\ell}{n_k}\nonumber\\[-8pt]\\[-8pt]
&\ge&
\liminf_{k\to\infty}\frac{n_k}{\tau_{n_k+1}}\frac{X_{\tau_{n_k}}
\cdot\ell}{n_k}
\ge\frac{1}{C},\nonumber
\end{eqnarray}
where $C$ is the constant appearing in Lemma~\ref{ltaun} and in
the last inequality we have used (\ref{xinf}) and
Lemma~\ref{ltaun}. This proves Theorem~\ref{tballisticity}.

\subsection{\texorpdfstring{Proof of Theorem \protect\ref{trwre}}{Proof of Theorem 1.2}}
As explained in the first paragraph of Section~\ref{srenewalstr},
let us note
that for each environment $\omega$, under the law $P_{\omega,0}$,
the ERWRE has a sequence of regeneration times $(\tau_n,n\ge1)$
which satisfy Propositions~\ref{tail-reg} and~\ref{conditioning}.

Furthermore (as explained, e.g., in~\cite{SZ99}), within
the context of random walk in random environment,
the following proposition is satisfied (the fact that
$P_0[\bar D_{0}=\infty]>0$ follows from Proposition~\ref{pgoaway}
below and the observation that a.s. the random walk satisfies
Conditions~\ref{ConditionB},~\ref{ConditionC+} and~\ref{ConditionE}):
\begin{prop}
\label{piid}
Let $(\tau_n,n\ge1)$ be the regeneration times of an ERWRE.
\begin{longlist}[(a)]
\item[(a)] Under\vspace*{1pt} the annealed law $P_0$,
$\tau_1,\tau_2-\tau_1,\tau_3-\tau_2,\ldots$ are independent
and $\tau_2-\tau_1,\tau_3-\tau_2,\ldots$ are i.i.d.
and with the same law as $\tau_1$ under
$P_0[ \cdot|\bar D_{0}=\infty]$.
\item[(b)] Under the annealed law $P_0$,
$X_{(\cdot\land\tau_1)}, X_{((\cdot+\tau_1)\land\tau_2)}-X_{\tau_1},
X_{((\cdot+\tau_2)\land\tau_3)}-X_{\tau_2},\ldots$ are independent,
and
$X_{((\cdot+\tau_1)\land\tau_2)}-X_{\tau_1},
X_{((\cdot+\tau_2)\land\tau_3)}-X_{\tau_2},\ldots$ are i.i.d.
with the same law as $X_{(\cdot\land\tau_1)}$ under
$P_0[ \cdot|\bar D_0=\infty]$.
\end{longlist}
\end{prop}

Propositions~\ref{tail-reg} and~\ref{piid} have the following two
corollaries.
\begin{cor} Let $A$ be a Borel subset of
$(\mathbb Z^d)^{\mathbb N}$.
Then for every $n\ge1$,
\[
P_0\bigl[X_{\tau_n+\cdot}-X_{\tau_n}\in A\mid\G^{(n)}\bigr]=
P_0[X_{\cdot}\in A\mid\bar D_0=\infty].\vadjust{\goodbreak}
\]
\end{cor}
\begin{cor}
\label{c2}
Let $(\tau_n,n\ge1)$ be the regeneration times
of the ERWRE. Then the following are satisfied:
\begin{longlist}[(a)]
\item[(a)] $E_0[\tau_1^2]<\infty$ and $E_0[\tau_1^2\mid\bar
D_0=\infty]<\infty$;
\item[(b)] $E_0[X^2_{\tau_1}]<\infty$ and
$E_0[X^2_{\tau_1}\mid\bar D_0=\infty]<\infty$.
\end{longlist}
\end{cor}

Now, using standard methods, one can prove
part (i) of Theorem~\ref{trwre}, showing that a.s.
\[
v=\lim_{n\to\infty}\frac{X_n}{n}=
\frac{E_0[X_{\tau_1}\mid\bar D_0=\infty]}{E_0[\tau_1\mid\bar
D_0=\infty]}.
\]
Indeed, by Proposition~\ref{piid} and Corollary~\ref{c2} we have
that
%
\begin{eqnarray}
\label{lim}
\lim_{n\to\infty}\frac{\tau_n}{n}&=&E_0[\tau_1\mid\bar D_0=\infty]
\quad\mbox{and}\nonumber\\[-8pt]\\[-8pt]
\lim_{n\to\infty}\frac{X_{\tau_n}}{n}&=&E_0[X_{\tau_1} \mid
\bar D_0=\infty].\nonumber
\end{eqnarray}
Then, standard arguments enable us to deduce part (i) of
Theorem~\ref{trwre} from (\ref{lim}).
The proof of part (ii) of Theorem~\ref{trwre} follows the
methods, for example, of~\cite{S00}, using Corollary~\ref{c2},
deducing that the covariance matrix of the limiting distribution is
given by
\[
A=\frac{E_0[(X_{\tau_1}-\tau_1v)^t(X_{\tau_1}-\tau_1v) \mid
\bar D_0=\infty]}{E_0[\tau_1 \mid\bar D_0=\infty]}.
\]

\section{Displacement estimates and tails of the regeneration times}
\label{srange}
In this section we will derive estimates
on the displacement of
generalized excited random walks with an arbitrary initial
cookie configuration. This will be used to prove
the tail estimates for the regeneration times in
Proposition~\ref{tail-reg}.
A~key ingredient in the proofs will be estimates on the range
of the process.

\subsection{Displacement estimates}
\label{sdisplest}
For the sake of completeness and for future reference, we introduce
also:
\renewcommand{\theCondition}{C}
\begin{Condition}\label{ConditionC}
Let $\ell\in\Sph^{d-1}$. We say that Condition~\ref{ConditionC} is satisfied with
respect to $\ell$ if
\[
\IE(X_{n+1}-X_n\mid\F_n)=0 \qquad\mbox{on }
\{\mbox{there exists }k<n \mbox{ such that }X_k=X_n\}
\]
and
\[
\IE(X_{n+1}-X_n\mid\F_n)\cdot\ell\geq0 \qquad\mbox{on }
\{X_k\neq X_n \mbox{ for all }k<n\}.
\]
\end{Condition}

That is, Condition~\ref{ConditionC} does not require that on the first visit to a
site the particle gets a \textit{uniformly} positive push in the
direction $\ell$.

Now, we need to consider processes starting at an arbitrary cookie
environment.\vadjust{\goodbreak}
\renewcommand{\theConditions}{C$_A$ \textup{and} C$^+_A$}
\begin{Conditions}\hypertarget{ConditionCA}
Let $\ell\in\Sph^{d-1}$
and $A\subset\Z^d$.
We say that Condition~\hyperlink{ConditionCA}{C$_A$} is satisfied with respect to $\ell$, if
\begin{eqnarray}
&&\IE(X_{n+1}-X_n\mid\F_n)=0 \hspace*{115pt}\nonumber\\
&&\eqntext{\mbox{on the event }
\{\mbox{there exists }k<n
\mbox{ such that }X_k=X_n \mbox{
or
}X_n\notin A\}}
\end{eqnarray}
and
\[
\IE(X_{n+1}-X_n\mid\F_n)\cdot\ell\geq0 \qquad\mbox{on }
\{X_k\neq X_n \mbox{ for all }k<n \mbox{ and }X_n\in A\}.
\]
If (in addition to the first display) there exist $\lambda> 0$ such
that
\[
\IE(X_{n+1}-X_n\mid\F_n)\cdot\ell\geq\lambda\qquad\mbox{on }
\{X_k\neq X_n \mbox{ for all }k<n \mbox{ and }X_n\in A\},
\]
we say that Condition \hyperlink{ConditionCA}{C$^+_A$} holds.
\end{Conditions}

Note that Condition \hyperlink{ConditionCA}{C$_A$} (\hyperlink{ConditionCA}{C$^+_A$})
becomes Condition~\ref{ConditionC} (\ref{ConditionC+}), if $A=\Z^d$. On the other hand,
Condition \hyperlink{ConditionCA}{C$_A$} with $A=\varnothing$ simply means that the process is
a $d$-dimensional martingale.
Throughout, given a stochastic process $({\tilde X}_n,n\ge0)$ on
the lattice $\Z^d$, we denote its range at time $n$ by
\[
\RR_n^{\tilde X}=\{x\in\Z^d\dvtx{\tilde X}_k=x \mbox{ for some }
0\le k\le n\}.
\]
Let us use also the notation $|U|$ for the cardinality of $U$,
where $U\subset\Z^d$.

Now, the key fact in this paper is that for \textit{any} process
that satisfies Conditions~\ref{ConditionB},~\ref{ConditionE} and
\hyperlink{ConditionCA}{C$_A$}
with an arbitrary $A\subset\Z^d$, the number of distinct sites
visited by time $n$ is typically much larger than $n^{1/2}$:
\begin{prop}
\label{phitpoints}
Suppose\vspace*{1pt} that Conditions~\ref{ConditionB},~\ref{ConditionE} and
\textup{\hyperlink{ConditionCA}{C$_A$}} (for some
$A\subset\Z^d$, $d\geq2$) hold for a process ${\tilde X}=({\tilde
X}_n,
n=0,1,2, \ldots)$.
Then, there exist positive constants $\alpha, \gamma_1,\gamma_2$
which depend only on $d,K,h,r$, such that
%
\begin{equation}
\label{ineqhitpoints}
\IP[|\RR_n^{\tilde X}|<n^{{1/2}+\alpha}] <
e^{-\gamma_1n^{\gamma_2}}
\end{equation}
for all $n\geq1$.
\end{prop}

The proof of this proposition is postponed
to Section~\ref{shitpoints}.

Now, let us recall Azuma's inequalities. Let $a>0$.
If $\{Z_n\}_{n\in\N}$ is a
martingale with respect to some filtration, and such that
$|Z_k-Z_{k-1}|<c$ a.s., then (cf., e.g., Theorem 7.2.1
of~\cite{AS00})
%
\begin{equation}
\label{Azuma}
\IP[|Z_n - Z_0| \geq a] \leq2\exp\biggl(-\frac{a^2}{2nc^2} \biggr).
\end{equation}
If $ \{\tilde Z_n\}_{n\in\N}$ is a super-martingale,
$|\tilde Z_k-\tilde Z_{k-1}|<c$ a.s., then (see, e.g., Lem\-ma~1
of~\cite{W95})
%
\begin{equation}
\label{Azuma1}
\IP[\tilde Z_n - \tilde Z_0 \geq a]
\leq\exp\biggl(-\frac{a^2}{2nc^2} \biggr).
\end{equation}
Let $H(a,b)\subset\Z^d $ be defined as
\[
H(a,b)=\{x\in\Z^d\dvtx x\cdot\ell\in[a,b]\}.\vadjust{\goodbreak}
\]
Now we obtain an important consequence of
Proposition~\ref{phitpoints}: if the cookies' configuration $A$ is
such that there are enough cookies around the starting point, then
the process is likely to advance in the direction $\ell$.\vspace*{-2pt}
\begin{prop}
\label{padvance}
Suppose that the process $X$ satisfies Conditions~\ref{ConditionB},~\ref{ConditionE},
\textup{\hyperlink{ConditionCA}{C$^+_A$}},
and there exists $\delta_0>0$ such that for some $n\geq1$
%
\begin{equation}
\label{condonA}
\bigl|(\Z^d\setminus A)\cap H\bigl( -n^{{1/2}+\delta_0},
\tfrac{2}{3}\lambda n^{{1/2}+\alpha}\bigr)\bigr|
\leq\tfrac{1}{3}n^{{1/2}+\alpha},
\end{equation}
where $\alpha$ is from Proposition~\ref{phitpoints}.
Then, for some positive constants $\gamma_3,\gamma_4$ that depend
only on $d,K,h,r,\lambda,\delta_0$, we have
%
\begin{equation}
\label{ineqadvance}
\IP\bigl[X_n\cdot\ell
< \tfrac{1}{3}\lambda n^{{1/2}+\alpha}\bigr]
< e^{-\gamma_3n^{\gamma_4}}.\vspace*{-2pt}
\end{equation}
\end{prop}
\begin{pf}
First, note that, by (\ref{Azuma1})
\[
\IP\biggl[\max_{k\leq n}X_k\cdot\ell> \frac{2}{3}\lambda
n^{{1/2}+\alpha}, X_n\cdot\ell< \frac{1}{3}\lambda
n^{{1/2}+\alpha}\biggr] \leq C_1ne^{-C_2n^{2\alpha}}
\]
for some $C_1,C_2>0$. Observe also that, again by (\ref{Azuma1}),
\[
\IP\Bigl[\min_{k\leq n}X_k\cdot\ell< -n^{{1/2}+\delta
_0}\Bigr]
\leq C_3ne^{-C_4n^{2\delta_0}}
\]
for some $C_3,C_4>0$.
Next, let
\[
\D_k = \IE(X_{k+1}-X_k\mid\F_k)
\]
be the (conditional) drift at time $k$. Then, the process
\[
Y_n = X_n - \sum_{k=0}^{n-1}\D_k
\]
is a martingale with bounded increments.
By (\ref{condonA}) and Condition \hyperlink{ConditionCA}{C$_A$}, on the event
\[
\{|\RR_n^X|\geq n^{{1/2}+\alpha}\} \cap
\bigl\{X_k\in H\bigl( -n^{{1/2}+\delta_0},
\tfrac{2}{3}\lambda n^{{1/2}+\alpha}\bigr)
\mbox{ for all } k\leq n\bigr\},
\]
we have
\[
\Biggl( \sum_{k=0}^{n-1}\D_k\Biggr)\cdot\ell>
\frac{2}{3}\lambda n^{{1/2}+\alpha}.
\]
Hence, using (\ref{Azuma}) and Proposition~\ref{phitpoints},
we conclude the proof of Proposition~\ref{padvance}.\vspace*{-2pt}
\end{pf}
\begin{rmk}
Condition (\ref{condonA}) is enough for our needs in this
paper, but, in fact, from the proof of Proposition~\ref{padvance}
one can see that it can be relaxed in the following way: if in the
strip of width $O(n^{{1/2}+\hat\alpha})$ (where
$0<\hat\alpha<\alpha$) there are not more than
$O(n^{{1/2}+\alpha})$ eaten cookies, then the\vspace*{1pt} process is
likely to advance by at least $O(n^{{1/2}+\hat\alpha})$ by
time $n$.\vspace*{-2pt}
\end{rmk}

For each $\ell\in\Sph^{d-1}$, define the half-space
$\MM_\ell=\{x\in\Z^d\dvtx x\cdot\ell>0\}$. Next, we prove that if
for some $\ell\in\Sph^{d-1}$ we know that\vadjust{\goodbreak}
all the cookies in the half-space~$\MM_\ell$ are present and that
Condition \hyperlink{ConditionCA}{C$_A$} is satisfied with respect to $\ell$, then with
uniformly positive probability the process never goes below its
initial location:\vspace*{-2pt}
\begin{prop}
\label{pgoaway}
Assume Conditions~\ref{ConditionB},~\ref{ConditionE}, \textup{\hyperlink{ConditionCA}{C$^+_A$}}, and suppose that $A$ is such
that $\MM_\ell\subset A$. Then there exists
$\psi=\psi(d,K,h,r,\lambda)>0$ such that
%
\begin{equation}
\label{ineqgoaway}
\IP[\bar D_0=\infty]\geq\IP[X_n\cdot\ell>0 \mbox{ for all
}n\geq1] \geq\psi.\vspace*{-2pt}
\end{equation}
%
\end{prop}
\begin{pf}
Clearly, since on $\{X_n\cdot\ell> 0 \mbox{ for all }n\geq1\}$
the process does not see the cookies' configuration on
$\Z^d\setminus\MM_\ell$, it is enough to prove this proposition for
the case $A=\Z^d$. For a (large enough) integer $m$,
consider the event
\[
U_0 = \{(X_{k+1}-X_k)\cdot\ell\geq r \mbox{ for all }
k=0,\ldots,\lceil r^{-1}\rceil m-1 \};
\]
observe that $X_{\lceil r^{-1}\rceil m}\cdot\ell\geq m$ on $U_0$.
By (\ref{snos}) we have for any fixed $m>0$
%
\begin{equation}
\label{otoshliot0}
\IP[U_0] \geq h^{\lceil r^{-1}\rceil m}.
\end{equation}

Let $A'=\Z^d\setminus\{X_0,\ldots,X_{\lceil r^{-1}\rceil m-1}\}$,
$y_0=X_{\lceil r^{-1}\rceil m}$, and abbreviate $W_k=X_{\lceil
r^{-1}\rceil m+k}$, $k\geq0$, so that $W_0=y_0$. Observe that,
if $m$ is large enough, then
the process $W-y_0$ satisfies Conditions~\ref{ConditionB} and~\ref{ConditionE}, and $A'-y_0$
satisfies (\ref{condonA}) for all $n\geq m^{2-\eps}$ for some
small enough $\delta_0$. Now, suppose
that the event $U_0$ occurs and let us fix $\eps$ such that
$(2-\eps)(\frac{1}{2}+\alpha)>1$. Denote $m_0=0, m_1=m$ and
$m_{k+1}=\frac{\lambda}{3}m_k^{(2-\eps)({1/2}+\alpha)}$
for $k\geq1$. Consider the events
\begin{eqnarray*}
G_k &=& \Bigl\{ \min_{j\leq m_k^{2-\eps}}(W_j-W_0)\cdot\ell
> -m_k\Bigr\},\\
U_k &=& \bigl\{W_{\lfloor m_k^{2-\eps}\rfloor} \geq
m_{k+1}\bigr\},\qquad
k\geq1.
\end{eqnarray*}
Observe that
%
\begin{equation}
\label{ubezhali}
\{X_n\cdot\ell>0 \mbox{ for all }n\geq1\} \supset
\Biggl(\bigcap_{k=1}^\infty(G_k\cap U_k)\Biggr)\cap U_0
\end{equation}
(indeed, on $G_k\cap U_{k-1}$ we have $X_n\cdot\ell>0$ for all
$n\in(m_{k-1}^{2-\eps},m_k^{2-\eps}]$). Also, by~(\ref{Azuma1}),
%
\begin{equation}
\label{ocGk}
\IP[G_k\mid U_0] \geq1-m_k^{2-\eps}e^{-2K^2m_k^\eps}
\end{equation}
and, by Proposition~\ref{padvance},
%
\begin{equation}
\label{ocUk}
\IP[U_k\mid U_0] \geq1-e^{-\gamma_3m_k^{(2-\eps)\gamma_4}}.
\end{equation}
Since $\IP[\bigcap_{k=1}^\infty(G_k\cap U_k)\mid U_0] \geq
1-\sum_{k=1}^\infty(\IP[G_k^c\mid U_0]+\IP[U_k^c\mid U_0])$,
Proposition~\ref{pgoaway} now follows from (\ref{otoshliot0})
and (\ref{ubezhali}).\vspace*{-2pt}
\end{pf}

\subsection{\texorpdfstring{Proof of Proposition \protect\ref{tail-reg}}{Proof of Proposition 2.1}}
\label{sprooftails}
Here we will prove Proposition~\ref{tail-reg}, closely following \cite
{BR07}. We will need the following result:\vadjust{\goodbreak}
\begin{prop}
\label{estimates}
Consider a generalized excited random walk with respect
to a filtration $(\F_k,k\ge0)$ and in the direction $\ell$.
Let $(\tau_k,k\ge1)$ be the corresponding regeneration times.
Then:
\begin{longlist}[(iii)]
\item[(i)] there exists a positive constant $\psi$ depending only
on $d,K,h,r,\lambda$, such that
%
\begin{equation}
\sup_{j,k\ge1} \IP\bigl[ D_j^{(k)}<\infty
\mid\G_j^{(k)}\bigr] < 1-\psi;\vspace*{-2pt}
\end{equation}
\item[(ii)] there exist positive constants $\gamma_3,\gamma_4$
depending only
on $d,K,h,r,\lambda$, such that
%
\begin{equation}
\label{ineqadvance2}
\sup_{k\ge1}\IP\biggl[(X_{\tau_k+n}-X_{\tau_k})\cdot\ell
< \frac{1}{3}\lambda n^{{1/2}+\alpha} \Bigm|
\G^{(k)}_0\biggr]
< e^{-\gamma_3n^{\gamma_4}};\vspace*{-2pt}
\end{equation}
\item[(iii)] there exist positive constants $\gamma_5,\gamma_6$
depending only
on $d,K,h,r,\lambda$, such that
%
\begin{equation}
\sup_{j,k\ge1} \IP\bigl[n\le
D_j^{(k)}-S_j^{(k)}<\infty\mid\G^{(k)}_0\bigr]
< e^{-\gamma_5n^{\gamma_6}}.\vspace*{-2pt}
\end{equation}
\end{longlist}
\end{prop}

Before proving Proposition~\ref{estimates},
let us see how it implies
Proposition~\ref{tail-reg}. We follow the proof of Proposition 1
in~\cite{BR07}. Let $a_1,a_2,a_3$ be positive real numbers
such that $a_1<\frac{1}{2}+\alpha$ and $a_2+a_3<a_1$.
For each integer $n>0$,
let $u_n=\lfloor n^{a_1}\rfloor$, $v_n=\lfloor n^{a_2}\rfloor$
and $w_n=\lfloor n^{a_3}\rfloor$.
We now choose $n$ large enough so that $(K+1)v_n(w_n+1)+2+K\le u_n$.
Let
\[
A_n=\{(X_{\tau_k+n}-X_{\tau_k})\cdot\ell\le u_n\},\qquad
B_n=\bigcap_{j=0}^{v_n}\bigl\{D_j^{(k)}<\infty\bigr\}\vspace*{-2pt}
\]
and
\[
F_n=\bigcup_{j=0}^{v_n}\bigl\{w_n\le D_j^{(k)}-S_j^{(k)}<\infty\bigr\}.\vspace*{-2pt}
\]
We will show that
%
\begin{equation}
\label{inclusion}
A_n^c\cap B_n^c\cap F_n^c\subset\{\tau_{k+1}-\tau_k<n\}.\vspace*{-2pt}
\end{equation}
To this end, for each natural $n\ge\tau_k$ define
\[
r_n=\max\{ (X_{j}-X_{\tau_k})\cdot\ell\dvtx\tau_k\le j\le n\}.\vspace*{-2pt}
\]
On the event $A_n^c\cap B_n^c\cap C_n^c$, we can define
\[
M=\inf\bigl\{0\le j\le v_n\dvtx D_j^{(k)}=\infty\bigr\},\vspace*{-2pt}
\]
and it is
true that $\tau_{k+1}=S_M^{(k)}$.
Hence, we need to prove that $S_M^{(k)}-\tau_k<n$.
Let us set $D^{(k)}_{-1}=\tau_k$.
By definition we know that $D^{(k)}_{M-1}<\infty$.
We now write
\[
r_{D^{(k)}_{M-1}}=\sum_{j=0}^{M-1}\bigl(r_{D^{(k)}_j}-
r_{S^{(k)}_j}\bigr)+\bigl(r_{S^{(k)}_j}-r_{D^{(k)}_{j-1}}\bigr)\vadjust{\goodbreak}
\]
with the
convention $\sum_{j=0}^{-1}=0$. Since, by Condition~\ref{ConditionB},
the range of each jump is at most $K$, we have that for each $0\le
j\le M-1$, $r_{D^{(k)}_j}-r_{S^{(k)}_j}\le K(D^{(k)}_j-S^{(k)}_j)$.
And by definition, for each $0\le j\le M-1$, we have that
$r_{S^{(k)}_j}-r_{D^{(k)}_{j-1}}\le K+1$.
But on the event $F_n^c$, since for each $0\le j\le M-1$ we have
$D^{(k)}_j<\infty$, it is true that $D^{(k)}_j-S^{(k)}_j\le w_n$.
It follows that
\[
r_{D^{(k)}_{M-1}}\le(K+1)v_n(w_n+1).\vspace*{-1pt}
\]
Since we have chosen $n$ sufficiently large so that
$(K+1)v_n(w_n+1)+2+K\le u_n$, we have
$r_{D^{(k)}_{M-1}}+2+K\le u_n$. Now, on $A^c_n$, we have
that $X_{\tau_k+n}\cdot\ell-X_{\tau_k}\cdot\ell>u_n$.
Hence, $X_{\tau_k+n}\cdot
\ell-X_{\tau_k}\cdot\ell>r_{D^{(k)}_{M-1}}+2+K$
and the smallest $i$ such that
$X_{\tau_k+i}\cdot\ell-X_{\tau_k}\cdot\ell
>r_{D^{(k)}_{M-1}}+1$ must be smaller than $n$. It follows
that\vspace*{0.5pt}
$S_M^{(k)}-\tau_k<n$, and this concludes the proof
of (\ref{inclusion}).

Now let us show that (\ref{inclusion}) is enough to prove
Proposition~\ref{tail-reg}.
Indeed, by parts (ii) and (iii) of Proposition~\ref{estimates},
we have that
%
\begin{equation}
\label{bound-an}
\IP\bigl[A_n \mid\G^{(k)}_0\bigr]\le e^{-\gamma_3n^{\gamma_4}}\vspace*{-1pt}
\end{equation}
and
%
\begin{equation}
\label{bound-cn}
\IP\bigl[F_n \mid\G^{(k)}_0\bigr]\le\frac{1}{\gamma_7}e^{-n^{\gamma_7}}\vspace*{-1pt}
\end{equation}
for some $\gamma_7$ such that $0<\gamma_7<\gamma_6$.
Furthermore, by part (i) of Proposition~\ref{estimates},
we see that
%
\begin{equation}
\label{bound-bn}
\IP\bigl[B_n \mid\G^{(k)}_0\bigr]\le(1-\psi)^{v_n}.\vspace*{-1pt}
\end{equation}
It is clear now that estimates (\ref{bound-an}),
(\ref{bound-cn}) and (\ref{bound-bn}) applied to
inclusion~(\ref{inclusion}), imply the statement of
Proposition~\ref{tail-reg}.
\begin{pf*}{Proof of Proposition~\ref{estimates}}

\textsc{Proof of part} (i).
Let $k$ and $j$ be fixed positive integers. By part (ii) of
Proposition~\ref{conditioning}, we have that
%
\begin{eqnarray}\label{pi}
&&\IP\bigl[D^{(k)}_{j}=\infty\mid{\G}^{(k)}_j\bigr]\nonumber\\[-1.5pt]
&&\qquad=\sum_{n=1}^\infty\sum_{m=1}^\infty{\mathbf{1}}_{\{\tau_k=n\}
}(\varpi)
{\mathbf{1}}_{\{S_j^{(k)}=n+m\}}(\varpi)\\[-1.5pt]
&&\qquad\quad\hphantom{\sum_{n=1}^\infty\sum_{m=1}^\infty}
{}\times\IP[\bar D_{n+m}=\infty\mid
\bar D_n=\infty,\F_{n+m}].\nonumber\vspace*{-1pt}
\end{eqnarray}
Now, on the paths such that
$\inf_{n\le j\le n+m}X_j\cdot\ell\ge X_n\cdot\ell$ (which happens when
$\tau_k=n$) we have that
\begin{eqnarray*}
\IP[\bar D_{n+m}=\infty\mid\bar D_n=\infty,\F_{n+m}]&=&
\frac{\IP[\bar D_{n+m}=\infty\mid\F_{n+m}]}
{\IP[\bar D_n=\infty\mid\F_{n+m}]}\\[-1.5pt]
&\ge&
\IP[\bar D_{n+m}=\infty\mid\F_{n+m}].\vadjust{\goodbreak}
\end{eqnarray*}
Hence, using Proposition~\ref{pgoaway} we get that
%
\begin{equation}\label{pii}
\IP[\bar D_{n+m}=\infty\mid\bar D_n=\infty, \F_{n+m}]
\ge\psi>0.
\end{equation}
Substituting (\ref{pii}) into (\ref{pi}) we obtain that
\[
\IP\bigl[D^{(k)}_{j}=\infty\mid{\G}^{(k)}_j\bigr]\ge\psi
\sum_{n=1}^\infty\sum_{m=1}^\infty{\mathbf{1}}_{\{\tau_k=n\}
}(\varpi)
{\mathbf{1}}_{\{S_j^{(k)}=n+m\}}(\varpi)\ge\psi>0.
\]

\textsc{Proof of part} (ii).
Note that
\begin{eqnarray*}
&&\IP\biggl[(X_{\tau_k+n}-X_{\tau_k})\cdot
\ell<\lambda\frac{1}{3}n^{{1/2}+\alpha}
\Bigm| \G^{(k)}_0 \biggr]\\
&&\quad=\sum_{m=1}^\infty
{\mathbf{1}}_{\{\tau_k=m\}}(\varpi)\IP\biggl[(X_{m+n}-X_m)\cdot
\ell<\lambda\frac{1}{3}n^{{1/2}+\alpha} \Bigm|
\bar D_m=\infty,\F_m\biggr]\\
&&\quad
\le
\sum_{m=1}^\infty
{\mathbf{1}}_{\{\tau_k=m\}}(\varpi)
\frac{1}{\IP[\bar D_m=\infty\mid\F_m]}
\IP\biggl[(X_{m+n}-X_m)\cdot
\ell<\lambda\frac{1}{3}n^{{1/2}+\alpha} \Bigm|
\F_m\biggr]\\
&&\quad\le
\sum_{m=1}^\infty
{\mathbf{1}}_{\{\tau_k=m\}}(\varpi)
\frac{1}{\IP[\bar D_m=\infty\mid\F_m]}
e^{-\gamma_3 n^{\gamma_4}}\\
&&\quad\le\frac{1}{\psi}
e^{-\gamma_3 n^{\gamma_4}},
\end{eqnarray*}
where in the first equality we have used part (i) of
Proposition~\ref{conditioning},
in the second to last inequality Proposition~\ref{padvance}
and in the last step Proposition~\ref{pgoaway}.

\textsc{Proof of part} (iii). This follows from part
(ii) in
analogy with the proof of Lemma 9 of~\cite{BR07}.
\end{pf*}
%
%
%

\section{On the number of distinct sites visited by the cookie
process}
\label{shitpoints}
In this section, after obtaining some auxiliary results, we prove
Proposition~\ref{phitpoints}.

\subsection{Some preliminary facts}
\label{sprelimfacts}
Denote by
\[
L_n(m)=\sum_{j=0}^{n}{\mathbf{1}}_{\{X_j\cdot\ell\in[m,m+1)\}}
\]
the local time on the corresponding strip. Let also
%
\begin{equation}
\label{deftauU}
\tau_U^X=\min\{n\ge0\dvtx X_n\in U\}
\end{equation}
be the entrance time to set $U$ for the process $X$.

In the sequel we will use the following simple facts: if $\tau$ is a
finite stopping time, and $X$ is a process that satisfies
Conditions~\ref{ConditionB},~\ref{ConditionE} and \hyperlink{ConditionCA}{C$_A$}, then the process $X\circ\theta_\tau$
also satisfies Conditions~\ref{ConditionB},~\ref{ConditionE} and \hyperlink{ConditionCA}{C$_A$} (but maybe with
different $A$, which is random and $\F_\tau$-measurable).
If $Y$ satisfies\vadjust{\goodbreak}
Conditions~\ref{ConditionB},~\ref{ConditionE} and \hyperlink{ConditionCA}{C$_A$} with $A=\varnothing$,
then $Y\circ\theta_\tau$
also satisfies Conditions~\ref{ConditionB},~\ref{ConditionE} and \hyperlink{ConditionCA}{C$_A$} with $A=\varnothing$.

We need the following.
\begin{lmm}
\label{llocaltime}
Under Conditions~\ref{ConditionB},~\ref{ConditionE}, \textup{\hyperlink{ConditionCA}{C$_A$}},
for any $\delta>0$ there exists a~constant $\gamma'_1$ such that for all $m$ we have
\[
\IP[L_n(m)\ge n^{{1/2}+2\delta}]\le
e^{-\gamma'_1n^\delta}.
\]
\end{lmm}
\begin{pf}
Without restricting generality,
suppose that $m=-1$.
Denote ${\hat T}_0=0$, and
\[
{\hat T}_{k+1}=\min\{j> {\hat T}_k\dvtx X_j\cdot\ell\in[-1,0)\}.
\]
With this notation,
%
\begin{equation}
\label{charactL}
L_n(-1)=\max\{k\dvtx {\hat T}_k\le n\}.
\end{equation}
By Condition~\ref{ConditionE}, there exist $C_1>0$ and $i_0\geq1$ such that
for any stopping time $\tau$
%
\begin{equation}
\label{advance*}
\IP[(X_{\tau+i_0}-X_\tau)\cdot\ell\ge2\mid\F_\tau]\ge C_1.
\end{equation}
Using the optional stopping theorem,
for any $m$ and for any $x$ such that $x\cdot\ell\ge1$, we obtain
that
%
\begin{eqnarray}
\label{ocenkaexcursion}
1&\leq&(n^{{1/2}+\delta}+K)\nonumber\\
&&{}\times
\IP\bigl[{\tau^{X}_{H(n^{{1/2}+\delta},
+\infty)}}\circ{\theta_{{\hat T}_m+i_0}<
\tau^{X}_{H(-\infty, 0)}}\circ\theta_{{\hat T}_m+i_0} \mid
\F_{{\hat T}_m+i_0},\\
&&\hspace*{193pt} X_{{\hat T}_m+i_0}\cdot\ell= x\bigr].\nonumber
\end{eqnarray}
Then, by (\ref{advance*}) and (\ref{ocenkaexcursion}),
for any $m$ it holds that
\[
\IP\bigl[X_{{\hat T}_m+k}\cdot\ell>0 \mbox{ for all } i_0\le k\le
\tau^{X}_{H(n^{{1/2}+\delta}, +\infty)}
\circ\theta_{{\hat T}_m}
\mid\F_{{\hat T}_m}\bigr]
\ge C_2 n^{-{1/2}-\delta}.
\]
That is, starting at any $y$ such that
$y\cdot\ell\in[-1,0)$, by (\ref{advance*}) with a uniformly positive
probability the particle advances by at least distance $2$ in
direction~$\ell$ during the first $i_0$ steps, so that it comes to
$H(1,+\infty)$; then, by (\ref{ocenkaexcursion}), the particle will
visit $H(n^{{1/2}+\delta}, +\infty)$ before\vspace*{1pt} coming back to the
``negative'' half-space with probability at least
$O(n^{-{1/2}-\delta})$. But, if the particle managed to visit
$H(n^{{1/2}+\delta}, +\infty)$, by (\ref{Azuma1}) it is quite
likely that it will take more than~$n$ time units to go back to
$H(-\infty, 0)$.
So, by (\ref{Azuma1}), we have
\[
\IP[{\hat T}_{k+i_0}-{\hat T}_k>n \mid\F_{{\hat T}_k}]\ge
C_2 n^{-{1/2}-\delta}
\bigl(1-e^{-n^{2\delta}/(2K^2)}\bigr)\geq C_3
n^{-{1/2}-\delta}.
\]
%
Thus, one can write
\[
\IP[\mbox{there exists }k\le n^{{1/2}+2\delta}-i_0
\mbox{ such that }{\hat T}_{k+i_0}-{\hat T}_k>n]\ge
1-e^{-C_4n^\delta},
\]
and so, using (\ref{charactL}), we prove Lemma~\ref{llocaltime}.
\end{pf}

Denote by $B(x,s)=\{y\in\Z^d\dvtx \|y-x\|\le s\}$ the discrete ball
centered in $x$ and with radius $s$.\vadjust{\goodbreak}

Consistently with the discussion in Section~\ref{sintro},
let $\{Y_n\}_{n\in N}$ be a process which satisfies
Conditions~\ref{ConditionB},~\ref{ConditionE}, \hyperlink{ConditionCA}{C$_A$} with $A=\varnothing$ (i.e., $Y$ is a
$d$-dimensional martingale with uniformly bounded increments and
some uniform ellipticity).

Now we obtain some properties of the process $Y$ needed in
the course of the proof of our results.
\begin{lmm}
\label{lsubmart}
There exist $b\in(0,1)$ and $\gamma'_2>0$
(depending only on $K$, $h$,~$r$) such that
%
\begin{equation}
\label{eqsubmart}
\IE(\|Y_{n+1}\|^b\mid\F_n) \geq
\|Y_n\|^b{\mathbf{1}}_{\{\|Y_n\|>\gamma'_2\}}.
\end{equation}
\end{lmm}
\begin{pf}
Observe that, for any $\delta>0$ there exists $\eps'>0$ such that
for all $u\in\R$ with the property $|u|\leq\eps'$ and all
$b\in(\frac{1}{2},1)$ we have
%
\begin{equation}
\label{Taylor}
(1+u)^{b/2} \geq
1+\frac{b}{2}u-(1+\delta)\frac{b}{4}\biggl(1-\frac{b}{2}\biggr)u^2.
\end{equation}
Now, abbreviate
\[
P_{y,y+z} = \IP[Y_{n+1}-Y_n = z
\mid\F_n, Y_n=y];
\]
take $y$ such that $\frac{K}{\|y\|}<\eps'$, and write
\begin{eqnarray*}
&&\IE(\|Y_{n+1}\|^b - \|Y_n\|^b
\mid\F_n, Y_n=y)
\\
&&\qquad= \sum_z P_{y,y+z}(\|y+z\|^b-\|y\|^b)\\
&&\qquad= \|y\|^b \sum_z
P_{y,y+z}\biggl(\biggl(\frac{\|y+z\|^2}{\|y\|^2}\biggr)^{b/2} - 1
\biggr)\\
&&\qquad= \|y\|^b \sum_zP_{y,y+z}\biggl(\biggl(1 +
\frac{2y\cdot z + \|z\|^2}{\|y\|^2}\biggr)^{b/2} - 1 \biggr).
\end{eqnarray*}
Denote by $\phi_{y,z}$ the angle between $y$ and $z$ regarded as
vectors in $\R^d$.
By~(\ref{Taylor}), Condition~\ref{ConditionB} and the fact that
$\sum_z zP_{y,y+z}=0$, we can write
\begin{eqnarray*}
&&\IE(\|Y_{n+1}\|^b - \|Y_n\|^b
\mid Y_n=y)\\
&&\qquad\geq\|y\|^b \sum_z P_{y,y+z}
\biggl(\frac{by\cdot z}{\|y\|^2} + \frac{b\|z\|^2}{2\|y\|^2} -
(1+\delta)b\biggl(1-\frac{b}{2}\biggr)
\frac{(y\cdot z)^2}{\|y\|^4} \\
&&\hphantom{\|y\|^b \sum_z P_{y,y+z}
\biggl(}
\qquad\quad{} - (1+\delta)b\biggl(1-\frac{b}{2}\biggr)
\frac{(y\cdot z)\|z\|^2}{\|y\|^4} \\
&&\hspace*{185.7pt}{} -
(1+\delta)\frac{b}{4}\biggl(1-\frac{b}{2}\biggr)
\frac{\|z\|^4}{\|y\|^4} \biggr)\\
&&\qquad\geq\frac{b}{\|y\|^{2-b}}\sum_z P_{y,y+z}\biggl(\frac{b}{2}\|z\|^2
- (1+\delta)\biggl(1-\frac{b}{2}\biggr)\|z\|^2\cos^2\phi_{y,z}\biggr)\\
&&\qquad\quad{} +C_1\frac{K^3}{\|y\|^{3-b}}
+C_2\frac{K^4}{\|y\|^{4-b}},
\end{eqnarray*}
where $C_1$ and $C_2$ are some
(not necessarily positive) constants.
Since $d\geq2$ and using uniform ellipticity, we obtain that
if $\delta$ is small enough, and $b$ is close enough to $1$, we have
for any $y$
\[
\sum_z P_{y,y+z}\biggl(\frac{b}{2}\|z\|^2
- (1+\delta)\biggl(1-\frac{b}{2}\biggr)\|z\|^2\cos^2\phi_{y,z}\biggr) >
\delta' > 0
\]
[use Condition~\ref{ConditionB} and (\ref{snoss}) with some $\ell'$ such that
$\ell'\cdot y=0$].
Together with the previous computation, this completes
the proof of Lemma~\ref{lsubmart}.
\end{pf}

Next, we prove the following fact ($b$ is from
Lemma~\ref{lsubmart}):
\begin{lmm}
\label{lgenlocaltime}
Assume that a process $Y$ satisfies Conditions~\ref{ConditionB},~\ref{ConditionE}, \textup{\hyperlink{ConditionCA}{C$_A$}} with
$A=\varnothing$, and suppose also that $Y_0=x_0$. Then,
for any $\delta>0$ there exists $\gamma'_3>0$ such that
for all $x_0,y_0\in\Z^d$ and for all $n$ we have
%
\begin{equation}
\label{eqgenlocaltime}
\IP\Biggl[\sum_{j=1}^n {\mathbf{1}}_{\{Y_j=y_0\}}>n^{
{b/2}+\delta}\Biggr]
\leq e^{-\gamma'_3n^\delta}.
\end{equation}
\end{lmm}
\begin{pf}
Without\vspace*{2pt} restriction of generality, we may assume that $x_0=y_0=0$.
Abbreviate\vspace*{2pt} ${\tilde\tau}=\tau^Y_{\Z^d\setminus
B(0,\gamma'_2+1)}$ and ${\tilde V}=\{Y_m\neq0
\mbox{ for all }1\leq m\leq
{\tilde\tau}\}$, where~$\gamma'_2$ is from Lemma~\ref{lsubmart}.
By uniform ellipticity, there exists $C_1>0$ such that
%
\begin{equation}
\label{ocenka1ststep}
\IP[{\tilde V}] > C_1.
\end{equation}
%
By the optional stopping theorem and Lemma~\ref{lsubmart} we have
\[
(\gamma'_2+1)^b\leq(\gamma'_2)^b + C_2 (n^{1/2}+K)^b
\IP\bigl[\tau^{Y}_{\Z^d\setminus
B(0,C_2n^{1/2})}\circ\theta_{{\tilde\tau}} <
\tau^{Y}_{B(0,\gamma'_2)}\circ\theta_{{\tilde\tau}}
\mid{\tilde V}, \F_{{\tilde\tau}}\bigr],
\]
where $C_2$ is a (large) constant to be chosen later.
This implies that
%
\begin{equation}
\label{ochitball}
\IP\bigl[\tau^{Y}_{\Z^d\setminus
B(0,C_2n^{1/2})}\circ\theta_{{\tilde\tau}} <
\tau^Y_{B(0,\gamma'_2)}\circ\theta_{{\tilde\tau}}
\mid{\tilde V}, \F_{{\tilde\tau}}\bigr] \geq
\frac{C_3}{n^{b/2}}
\end{equation}
for some positive $C_3$ depending on $C_2$. Next we assume
that $C_2$ is sufficiently large
so that (\ref{Azuma}) implies that
for any stopping time ${\hat\tau}$
%
\begin{eqnarray}
\label{eqgoback}
&&\IP[\tau^{Y}_0\circ\theta_{\hat\tau}>n \mid\F_{\hat\tau},
Y_{\hat\tau}=y]\nonumber\hspace*{-35pt}\\[-8pt]\\[-8pt]
&&\qquad\ge 1-2\exp\biggl(-\frac{(C_2n^{1/2})^2}{2nK^2}\biggr)
\geq\frac{1}{2}\qquad
\mbox{for any }
y\in\Z^d\setminus B(0,C_2n^{1/2})\nonumber\hspace*{-35pt}
\end{eqnarray}
(to reach $0$ from $y$, the martingale $Y$ has to advance by at
least $C_2n^{1/2}$ units in some fixed direction).
Now, (\ref{ocenka1ststep}), (\ref{ochitball})
and (\ref{eqgoback}) imply that
%
\begin{equation}
\label{excursionZd}
\IP[Y_m\neq0 \mbox{ for all }m=1,\ldots,n] \geq
\frac{C_1C_3}{2n^{b/2}}.
\end{equation}
Then, proceeding as in the proof of Lemma~\ref{llocaltime}
[the argument after (\ref{ocenkaexcursion}) up to the end of the
proof], we obtain that (\ref{excursionZd})
implies (\ref{eqgenlocaltime}).
\end{pf}

Next, we prove that the process $Y$
typically hits sets which
contain enough points close to the starting place of the process:
\begin{lmm}
\label{lhitset}
Assume that a process $Y$ satisfies Conditions~\ref{ConditionB},~\ref{ConditionE},
\textup{\hyperlink{ConditionCA}{C$_A$}} with
\mbox{$A=\varnothing$}, and suppose that $Y_0=x$.
Consider an arbitrary $\delta>0$ and a~set~$U$ and
suppose that $|B(x,m^{1/2})\setminus U|\le
m^{1-{b/2}-2\delta}$, for some $m$. Then,
there exists $\gamma'_4>0$ such that
\[
\IP[\tau^Y_U\ge m^{1-\delta}]\le e^{-\gamma'_4m^\delta}.
\]
\end{lmm}
\begin{pf}
First, by (\ref{Azuma}), we have that
\[
\IP[Y_k\in B(x,m^{1/2}) \mbox{ for all }k\leq
m^{1-\delta}] \geq1-e^{-Cm^\delta}.
\]
Then, by Lemma~\ref{lgenlocaltime}, with probability at least
$1-e^{-Cm^\delta}$ by time $m^{1-\delta}$ every site
from $B(x,m^{1/2})$ will be visited less than
$m^{{b/2}+\delta}$ times, so we have
%
\[
|\RR^Y_{m^{1-\delta}}| = |\{Y_0,\ldots,Y_{m^{1-\delta}}\}|
> \frac{m^{1-\delta}}{m^{{b}/{2}+\delta}} =
m^{1-{b/2}-2\delta}
\]
with probability\vspace*{1pt} at least $1-e^{-C'm^\delta}$.
To complete the proof of Lemma~\ref{lhitset}, it remains to
observe that, since $|B(x,m^{1/2})\setminus U|\le
m^{1-{b/2}-2\delta}$, on the event
\[
\{|B(x,m^{1/2})\cap\RR^Y_{m^{1-\delta}}| >
m^{1-{b/2}-2\delta}\}
\]
we have $\{Y_0,\ldots,Y_{m^{1-\delta}}\}\cap U
\neq\varnothing$.
\end{pf}

\subsection{\texorpdfstring{Proof of Proposition \protect\ref{phitpoints}}{Proof of Proposition 4.1}}
\label{sproofhitpoints}
Fix $a\in(0,\frac{1}{2})$ and $\eps>0$ in such a way that
$(1-a+\eps)\wedge(\frac{1}{2}+\frac{a}{2}(1-b)-4\eps)>\frac{1}{2}$,
where $b$ is from Lemma~\ref{lsubmart};
also,\vspace*{1pt} fix $n\geq1$. In the rest of this section, we will not
explicitly indicate the dependence on $a,b,\eps,n$ by
sub/superscripts in the notation.
Let us denote by
\[
H_j=H\bigl(2(j-1)n^{{a/2}}, 2(j+1)n^{{a/2}}\bigr)
\]
the\vspace*{1pt} corresponding strip of width $4n^{{a/2}}$. We
say that the strip $H_j$ is a \textit{trap}
if $|\RR^X_n\cap H_j|\geq n^{a(1-{b/2})-2\eps}$.

Consider the event
\[
G = \bigl\{|\RR^X_n| \geq\tfrac{1}{2}
n^{(1-a+\eps)\wedge({1/2}+({a/2})(1-b)-4\eps)}\bigr\}.
\]
We are going to prove that
%
\begin{equation}
\label{manycookies}
\IP[G] \geq1 - e^{-C_1 n^{\eps/2}},\vadjust{\goodbreak}
\end{equation}
thus establishing Proposition~\ref{phitpoints}.
Let us introduce the event
\[
G_1 = \{L_n(k) \leq n^{{1/2}+\eps} \mbox{ for all }
k\in[- Kn, Kn]\}.
\]
By Lemma~\ref{llocaltime}, it holds that
%
\begin{equation}
\label{ocenkalocaltimes}
\IP[G_1] \geq1 - (2Kn+1)e^{-\gamma'_1 n^{\eps/2}}.
\end{equation}
Next, denote $\sigma_0=0$, and, inductively (assuming, of course,
that
$\lfloor n^{a-\eps}\rfloor\geq1$),
%
\begin{equation}
\label{defsigmak}\qquad
\sigma_{k+1} = \min\bigl\{j\geq\sigma_k+
\lfloor n^{a-\eps}\rfloor\dvtx |\RR^X_j\cap
B(X_j,n^{a/2})|\leq n^{a(1-{b/2})-2\eps}\bigr\}
\end{equation}
(formally, if such $j$ does not exist, we put $\sigma_{k+1}=\infty$).
Consider the event (to hit a new point means to visit a previously
unvisited site)
\begin{eqnarray*}
G_2&=&\bigl\{\mbox{at least one new point is hit on each
}\\
&&\hphantom{\bigl\{}
\mbox{of the time intervals }
[\sigma_{j-1},\sigma_j), j=1,\ldots,\tfrac{1}{2}n^{1-a+\eps}\bigr\}.
\end{eqnarray*}
Now, the key observation is the following: when the
process is walking on previously visited sites, it has zero drift.
So, if we only want to assure that at least one new
point it visited, this is equivalent to
considering the first moment when the process $Y$ (the
process without cookies) enters the previously unvisited set.
Then, by Lemma~\ref{lhitset} we have
\begin{eqnarray*}
&&\IP\bigl[\mbox{at least one new point is hit on each
} \\
&&\hphantom{\IP\bigl[}\mbox{of the time intervals }
[\sigma_{j-1},\sigma_j), j=1,\ldots,k\bigr] \\
&&\qquad\geq
1-ke^{-\gamma'_4n^{\eps/a}},
\end{eqnarray*}
so we obtain that
%
\begin{equation}
\label{ocenkanumbcookies}
\IP[G_2] \geq1 - \tfrac{1}{2}n^{1-a+\eps} e^{-\gamma'_4n^{\eps/a}}.
\end{equation}

Next, assuming that $n$ is so large that $8n^{1-\eps}<n/2$,
let us show that $(G_1\cap G_2)\subset G$.
Indeed, suppose that both $G_1$ and $G_2$ occur, but
$|\RR_n^X|<n^{{1/2}+{a/2}(1-b)-4\eps}$. Denote by
\[
{\hat L}_j = \sum_{k=2(j-1)n^{a/2}}^{2(j+1)n^{a/2}-1} L_n(k)
\]
the total number of visits to $H_j$. Then, on
$\{|\RR^X_n|<n^{{1/2}+{a/2}(1-b)-4\eps}\}$ the number of
traps
is at most $2n^{{1/2}-{a/2}-2\eps}$.
On the event $G_1$,
we can write
\[
\sum_j {\hat L}_j{\mathbf{1}}_{\{H_j\ \mathrm{is}\ \mathrm{a}\ \mathrm{trap}\}}\leq
4n^{{a/2}}\times2 n^{{1/2}-{a/2}-2\eps}\times
n^{{1/2}+\eps} = 8n^{1-\eps}.
\]
On the other hand, note that, since for $j\leq n$ we have
$\RR_j^X\subset\RR_n^X$, if
$|\RR_j^X\cap B(X_j,n^{a/2})| > n^{a(1-{b/2})-2\eps}$
then\vspace*{1pt} $X_j$
must be in a\vadjust{\goodbreak}
trap. Since $n$ is such that $8n^{1-\eps}<n/2$, we
obtain that, on the event
\[
\biggl\{\sum_j {\hat L}_j{\mathbf{1}}_{\{H_j\ \mathrm{is}\ \mathrm{a}\ \mathrm{trap}\}} \leq
8n^{1-\eps}\biggr\}
\]
we have that $\sigma_{n^{1-a+\eps}/2}<n$
(indeed, the total time spent in nontraps is at least $n/2$;
on the other hand, from the definition (\ref{defsigmak})
one can see that up to the moment $\sigma_k$ we can have at most
$kn^{a-\eps}$ instances $j$ such that $|\RR^X_j\cap
B(X_j,n^{a/2})|\leq n^{a(1-{b/2})-2\eps}$). But then,
on the event $G_2$ we have that $|\RR^X_n|\geq\frac{1}{2}
n^{(1-a+\eps)}$.
So,\vspace*{1pt} indeed $(G_1\cap G_2)\subset G$, and (\ref{manycookies})
follows from (\ref{ocenkalocaltimes})
and (\ref{ocenkanumbcookies}).
The proof of Proposition~\ref{phitpoints} is finished.

\section*{Acknowledgment}
The authors thank the referee for his careful
reading of the paper and valuable comments and suggestions.



%
\printaddresses

\end{document}